\RequirePackage{ifpdf}
\ifpdf 
\documentclass[pdftex]{sigma}
\else
\documentclass{sigma}
\fi

\newcommand{\res}{\operatorname{Res}}

\begin{document}

\allowdisplaybreaks

\renewcommand{\PaperNumber}{115}

\FirstPageHeading

\ShortArticleName{A Connection Formula of the Hahn--Exton $q$-Bessel Function}

\ArticleName{A Connection Formula\\ of the Hahn--Exton $\boldsymbol{q}$-Bessel Function}

\Author{Takeshi MORITA}

\AuthorNameForHeading{T.~Morita}

\Address{Graduate School of Information Science and Technology, Osaka University,\\
1-1  Machikaneyama-machi, Toyonaka, 560-0043, Japan}
\Email{\href{mailto:t-morita@cr.math.sci.osaka-u.ac.jp}{t-morita@cr.math.sci.osaka-u.ac.jp}}

\ArticleDates{Received May 11, 2011, in f\/inal form December 14, 2011; Published online December 16, 2011}

\Abstract{We show a connection formula of the Hahn--Exton $q$-Bessel function around the origin and the inf\/inity. We introduce the $q$-Borel transformation and the $q$-Laplace transformation following C.~Zhang to obtain the connection formula. We consider the limit $p\to 1^-$ of the connection formula.}

\Keywords{Hahn--Exton $q$-Bessel function; $q$-Borel transformation; connection problems}

\Classification{33D15; 34M40; 39A13}

\section{Introduction}
In this paper, we show a connection formula of the Hahn--Exton $q$-Bessel function $J_\nu^{(3)}(x;q)$. At f\/irst, we review the Bessel function and $q$-analogues of the Bessel function. The Bessel equation
\begin{gather*}
\frac{d^2u}{dz^2}+\frac{1}{z}\frac{du}{dz}+\left(1-\frac{\nu^2}{z^2}\right)u=0
\end{gather*}
has a solution $u(z)=J_\nu (z)$, $J_{-\nu}(z)$. Here, the Bessel function $J_\nu (z)$ is
\[J_\nu (z)=\frac{1}{\Gamma (\nu +1)}\left(\frac{z}{2}\right)^\nu{}_0F_1\left(-,\nu +1, -\frac{z^2}{4}\right).\]
The degenerated conf\/luent hypergeometric function ${}_0F_1(-,\alpha , z)$ is def\/ined by
\[{}_0F_1(-,\alpha , z)=\sum_{n\ge 0}\frac{1}{(\alpha )_n n!}z^n,\qquad (\alpha )_n=\alpha \{\alpha +1\}\cdots \{\alpha +(n-1)\}.\]
Both $J_\nu(z)$ and $J_{-\nu}(z)$ are linearly independent if $\nu\not\in\mathbb{Z}$.

It is known that there exists three dif\/ferent $q$-analogues of the Bessel function.
\begin{gather*}
J_\nu^{(1)}(x;q) :=\frac{(q^{\nu +1};q)_\infty}{(q;q)_\infty}\left(\frac{x}{2}\right)^\nu {}_2\varphi_1\left(0,0;q^{\nu +1};q,-\frac{x^2}{4}\right),\qquad |x|<2,\\
J_\nu^{(2)}(x;q) :=\frac{(q^{\nu +1};q)_\infty}{(q;q)_\infty}\left(\frac{x}{2}\right)^\nu {}_0\varphi_1\left(-;q^{\nu +1};q,-\frac{q^{\nu -1}x^2}{4}\right),\qquad x\in\mathbb{C},\\
J_\nu^{(3)}(x;q) :=\frac{(q^{\nu +1};q)_\infty}{(q;q)_\infty}x^\nu {}_1\varphi_1\left(0;q^{\nu +1};q,qx^2\right),\qquad x\in\mathbb{C}.
\end{gather*}
Here,
\begin{gather*}
(a;q)_n:=
\begin{cases}
1,&\quad n=0,\\
(1-a)(1-aq)\cdots (1-aq^{n-1}),&\quad n\ge 1,
\end{cases}
\\
(a;q)_\infty =\lim_{n\to\infty}(a;q)_n
\end{gather*}
and
\[(a_1,a_2,\dots ,a_m;q)_\infty =(a_1;q)_\infty (a_2;q)_\infty\cdots (a_m;q)_\infty .\]
Moreover, the basic hypergeometric series ${}_r\varphi_s$ is
\[{_r\varphi_s}(a_1,\dots ,a_r;b_1,\dots ,b_s;q,x)
:=\sum_{n\ge 0}\frac{(a_1,\dots ,a_r;q)_n}{(b_1,\dots ,b_s;q)_n(q;q)_n}
\left[(-1)^nq^{\frac{n(n-1)}{2}}\right]^{1+s-r}x^n.\]

The f\/irst and the second one
are called Jackson's f\/irst and second $q$-Bessel function and the third one is called the Hahn--Exton $q$-Bessel function.
They satisfy the following $q$-dif\/ference equations:
\begin{gather}
J^{(1)}_\nu :\quad u(xq) -\big(q^{\nu/2} +q^{-\nu/2}\big) u(xq^{1/2}) +\left(1+\frac{x^2}4\right)u(x)=0,\notag 
\\
J^{(2)}_\nu :\quad \left(1+\frac{qx^2}4\right)u(xq) -\big(q^{\nu/2} +q^{-\nu/2}\big) u\big(xq^{1/2}\big) +u(x)=0,\notag \\
J^{(3)}_\nu :\quad u(xq) -\left\{(q^{\nu/2} +q^{-\nu/2}) -q^{-\nu/2+1}x^2\right\} u\big(xq^{1/2}\big) + u(x)=0. \label{j3}
\end{gather}
The limits of these $q$-analogues of the Bessel function are the Bessel function when $q\to 1^-$:
\[\lim_{q\to 1^-}J_\nu^{(k)}\left((1-q)x;q\right)=J_{\nu}(x),\qquad k=1,2\]
and
\[\lim_{q\to 1^-}J_\nu^{(3)}\left((1-q)x;q\right)=J_\nu (2x).\]
The relation between $J_\nu^{(1)}(x;q)$ and $J_\nu^{(2)}(x;q)$ was found by Hahn \cite{Hahn} as follows:
\begin{gather}\label{J1J2}
J_\nu^{(2)}(x;q)=\left(-\frac{x^2}{4};q\right)_\infty J_\nu^{(1)}(x;q).
\end{gather}

Connection problems of the $q$-dif\/ference equation between the origin and the inf\/inity are studied by G.D.~Birkhof\/f \cite{Birkhoff}. We review connection formulae for several $q$-dif\/ference functions.
\begin{enumerate}\itemsep=0pt
\item Watson's formula.
In 1910 \cite{W}, Watson showed the connection formula of the basic hypergeometric function~${}_2\varphi_1$ as follows:
\begin{gather*}
{}_2 \varphi_1\left(a,b;c;q;x \right)=
\frac{(b,c/a;q)_\infty (a x,q/ a x;q)_\infty }{(c, b/a;q)_\infty (  x,q/   x;q)_\infty }
{}_2 \varphi_1\left(a,aq/c;aq/b;q;cq/abx \right) \nonumber \\
\hphantom{{}_2 \varphi_1\left(a,b;c;q;x \right)=}{} +
\frac{(a,c/b;q)_\infty (b x, q/ b x;q)_\infty }{(c, a/b;q)_\infty (  x,q/   x;q)_\infty }
{}_2 \varphi_1\left(b, bq/c; bq/a; q; cq/abx \right).
\end{gather*}
\item Connection formula of $J_\nu^{(1)}(x;q)$.
C.~Zhang has given some connection formulae for the solutions of the $q$-dif\/ference equations of conf\/luent type~\cite{Z0,Z1} and~\cite{Z2}. In \cite{Z1}, Zhang
has shown connection formulae for $J_\nu^{(1)}(x;q)$ and $J_\nu^{(2)}(x;q)$. The connection formula of $J_\nu^{(1)}(x;q)$ is given by
\begin{gather}
\frac{\left(\frac{\alpha}{\sqrt{p}x};p\right)_\infty}{\theta_p\left(-\frac{\alpha}{x}\right)}{}_2\varphi_1 \left(p^{\nu +\frac{1}{2}},p^{-\nu +\frac{1}{2}};-p;p,\frac{\alpha}{\sqrt{p}x}\right)\notag\\
\qquad{} =\frac{1}{\theta_p\left(-\frac{\alpha}{x}\right)}
 \left\{\frac{\theta_p\left(-\frac{\alpha q^{\frac{\nu}{2}}}{x}\right)}{(q,q^{-\nu};q)_\infty}{}_2\varphi_1\left(0,0;q^{\nu +1};q,-\frac{x^2}{4}\right)\right.\notag\\
\left.\qquad\quad{} + \frac{\theta_p\left(-\frac{\alpha q^{-\frac{\nu}{2}}}{x}\right)}{(q,q^{\nu};q)_\infty}{}_2\varphi_1\left(0,0;q^{-\nu +1};q,-\frac{x^2}{4}\right)
\right\},\label{Zhb}
\end{gather}
where $q=p^2$ and $\alpha^2=-4q^{3/2}$.
\end{enumerate}
The connection formula of $J_\nu^{(2)}(x;q)$ is obtained by (\ref{Zhb}) and (\ref{J1J2}).
But it is not known the connection formula of the Hahn--Exton $q$-Bessel function.

The Hahn--Exton $q$-Bessel equations (\ref{j3}) has two analytic solutions $u(x)=J_\nu^{(3)}(x)$, $J_{-\nu}^{(3)}(xp^{-\nu})$ around $x=0$ and has one analytic solution
$z\left(1/x\right)
=\frac{1}{\theta_p\left(-p^{\nu +2}/x\right)}\sum\limits_{n\ge 0}a_nx^{-n}$, $a_0=1$. We show a~connection formula of $J_\nu^{(3)}(x;q)$ in Section~\ref{section2} as follows:
\begin{theorem}\label{theorem1}
For any
$x\in\mathbb{C}^*\setminus [p^{\nu +2};p]$,
\begin{gather}
z\left(\frac{1}{x}\right) =\frac{1}{(p^{-2\nu},p;p)_\infty}\frac{\theta
_p\left(-\frac{p^{2\nu +2}}{x}\right)}{
\theta_p\left(-\frac{p^{\nu +2}}{x}\right)}\
_1\varphi_1\left(0,p^{1+2\nu};p,x\right)\notag\\
\phantom{z\left(\frac{1}{x}\right) =}{} +\frac{1}{(p^{2\nu},p;p)_\infty}\frac{\theta _p\left(-\frac{p^2}{x}\right)}{
\theta_p\left(-\frac{p^{\nu +2}}{x}\right)}\
_1\varphi_1\left(0,p^{1-2\nu};p,p^{-2\nu}x\right).\label{thm1}
\end{gather}
\end{theorem}
Here, $\theta_p(\cdot )$ is the theta function of Jacobi and $[\lambda ;q]$ is the $q$-spiral (see Section~\ref{section2}).
We use the $q$-Borel transformation and the $q$-Laplace transformation which is def\/ined by C.~Zhang in~\cite{Z1}.

In Section~\ref{section3}, we consider the limit $p\to 1^-$ of the connection formula. If we take a~suitable limit $p\to 1^-$ of (\ref{thm1}), we obtain
\[
H_\nu^{(2)}\left(\sqrt{z}\right)=\frac{-ie^{\nu\pi i}}{\sin\nu\pi}\left\{J_\nu\left(\sqrt{z}\right)-e^{-\nu\pi i}J_{-\nu}\left(\sqrt{z}\right)\right\}.
\]
Here, $H_\nu^{(2)}(z)$ is the Hankel function of the second kind. Thus we obtain a connection formula of the Bessel function as a limit $p\to 1^-$ of~(\ref{thm1}).

\section{The connection formula}\label{section2}

In this section, we give a connection formula of the Hahn--Exton $q$-Bessel function. We introduce the $p$-Borel transformation and the $p$-Laplace transformation to obtain the connection formula between the origin and the inf\/inity. These transformations are useful to consider connection problems. We assume that $q\in\mathbb{C}^*$ satisf\/ies $0< |q|<1$ and $q=p^2$. The $q$-dif\/ference opera\-tor~$\sigma_q$ is given by $\sigma_qf(x)=f(qx)$.

\subsection{The theta function of Jacobi}\label{section2.1}
Before we study connection problems, we review the theta function of Jacobi. The theta function of Jacobi is given by the following series:
\begin{definition}For any $x\in\mathbb{C}^*$,
\[
\theta_q(x)=\theta (x):=\sum_{n\in\mathbb{Z}}q^{\frac{n(n-1)}{2}}x^n.
\]
\end{definition}
We denote by $\theta_q(x)$ or more shortly $\theta (x)$. The theta function satisf\/ies Jacobi's triple product identity:
\[
\theta (x)=\left(q,-x,-\frac{q}{x};q\right)_\infty .
\]
The theta function satisf\/ies the $q$-dif\/ference equation as follows
\[
\theta (q^kx)=q^{-\frac{k(k-1)}{2}}x^{-k}\theta (x),\qquad \forall \, x\in\mathbb{C}^*.
\]
The theta function has the inversion formula
$x\theta(1/x)=\theta (x)$. For all f\/ixed $\lambda\in\mathbb{C}^*$, we def\/ine a~$q$-spiral $[\lambda ;q]:=\lambda q^{\mathbb{Z}}=\{\lambda q^k:k\in\mathbb{Z}\}$.
We remark that $\theta \left(\lambda q^k/x\right)=0$ if and only if $x\in [-\lambda ;q]$.

\subsection[The Hahn-Exton $q$-Bessel function]{The Hahn--Exton $\boldsymbol{q}$-Bessel function}\label{section2.2}

The Hahn--Exton $q$-Bessel function is def\/ined by
\[
J_\nu^{(3)}(x;q):=\frac{(q^{\nu +1};q)_\infty}{(q;q)_\infty}x^\nu\sum_{n\ge 0}\frac{(-1)^nq^{\frac{n(n-1)}{2}}}{(q^{\nu +1},q;q)_n}\left(qx^2\right)^n.
\]
The function $J_\nu^{(3)}(x;q)$ satisf\/ies the $q$-dif\/ference equation
\begin{gather}\label{heqb}
\left[\sigma_p^2- \left\{(p^\nu +p^{-\nu })-x^2p^{2-\nu}\right\}\sigma_p +1\right]y(x)=0.
\end{gather}
If we replace $\nu$ by $-\nu$ and $x$ by $xp^{-\nu}$,
we obtain $J_{-\nu}^{(3)}(xp^{-\nu};q)$ which is another solution of~(\ref{heqb}) around the origin. This solution corresponds to the classical Neumann function $Y_\nu (x)$ \cite{Z3}. We consider the behavior of equation~(\ref{heqb}) around the inf\/inity. We set $1/t$, formally $t^2\mapsto t$ and $z(t)=y(1/t)$. Then $z(t)$ satisf\/ies
\begin{gather}
\left[\sigma_p^2-\left\{(p^\nu +p^{-\nu})-\frac{p^{-2-\nu }}{t}\right\}\sigma_p+1\right]z\left(t\right)=0.
\label{b3}
\end{gather}

We set $\mathcal{E}(t)=1/\theta_p(-p^{\nu +2}t)$ and $f(t)=\sum\limits_{n\ge 0}a_nt^n$, $a_0=1$. We assume that $z(t)$ can be described as
\[z(t)=\mathcal{E}(t)f(t)=\frac{1}{\theta_p(-p^{\nu +2}t)}\left(\sum_{n\ge 0}a_nt^n\right).\]
Since $\mathcal{E}(t)$ satisf\/ies the following $q$-dif\/ference equation
\[\sigma_p\mathcal{E}(t)=-p^{\nu +2}t\mathcal{E}(t),\qquad
\sigma_p^2\mathcal{E}(t)=p^{2\nu +5}t^2\mathcal{E}(t),\]
we can check out that the function $f(t)$ satisf\/ies the equation
\begin{gather}
\left\{p^{2\nu +5}t^2\sigma_p^2 +p^{\nu +2} (p^\nu +p^{-\nu})t\sigma_p-\sigma_p +1\right\}f(t)=0.
\label{b4}
\end{gather}

\subsection[The $p$-Borel transformation and the $p$-Laplace transformation]{The $\boldsymbol{p}$-Borel transformation and the $\boldsymbol{p}$-Laplace transformation}\label{section2.3}

We def\/ine the $p$-Borel transformation and the $p$-Laplace transformation to solve the equation~(\ref{b4}), following~Zhang~\cite{Z1}.

\begin{definition}For $f(t)=\sum\limits_{n\ge 0}a_nt^n$, the $p$-Borel transformation is def\/ined by
\[g(\tau )=\left(\mathcal{B}_pf\right)(\tau ):=\sum_{n\ge 0}a_np^{-\frac{n(n-1)}{2}}\tau^n,\]
and the $p$-Laplace transformation is given by
\[
\left(\mathcal{L}_pg\right)(t):=\frac{1}{2\pi i}\int_{|\tau |=r}g(\tau )\theta_p\left(\frac{t}{\tau}\right)\frac{d\tau}{\tau }.
\]
Here, $r_0>0$ is enough small number.
\end{definition}
 The $p$-Borel transformation is considered as a formal inverse of the $p$-Laplace transformation.
\begin{lemma}We assume that the function $f$ can be $p$-Borel transformed to the analytic func\-tion~$g(\tau )$ around $\tau =0$. Then,
\[
\mathcal{L}_p\circ\mathcal{B}_pf=f.
\]
\end{lemma}

\begin{proof}We can prove this lemma calculating residues of the $p$-Laplace transformation around the origin.
\end{proof}

The $p$-Borel transformation has the following operational relation.
\begin{lemma}\label{oprel}For any $l,m \in \mathbb{Z}_{\ge 0}$,
\[
\mathcal{B}_p\big(t^m\sigma_p^l\big)=p^{-\frac{m(m-1)}{2}}\tau^m\sigma_p^{l-m}\mathcal{B}_p.
\]
\end{lemma}
Applying the $p$-Borel transformation to the equation (\ref{b4}) and using Lemma~\ref{oprel}, $g(\tau )$ satisf\/ies the f\/irst order dif\/ference equation
\[
g(p\tau )=\left(1+p^{2\nu +2}\tau \right)\left(1+p^2\tau \right)g(\tau ).
\]
Since $g(0)=1$, we get an inf\/inite product of $g(\tau )$:
\[
g(\tau )=\frac{1}{(-p^{2\nu +2}\tau ;p)_\infty(-p^2\tau ;p)_\infty}.
\]
Then $g(\tau )$ has single poles at
\[
\left\{-p^{-2\nu -2-k},-p^{-2-k};k\in\mathbb{Z}_{\ge 0}\right\}.
\]
We set
\[0<r<r_0:=\min\left\{\frac{1}{|p^{2\nu +2}|},  \frac{1}{|p^2|}\right\}.\]
and choose the radius $r>0$ such that $0<r<r_0$.
By Cauchy's residue theorem, the $p$-Laplace transform of $g(\tau )$ is
\begin{gather*}
f(t)= \frac{1}{2\pi i}\int_{|\tau |=r}g(\tau )\theta_p\left(\frac{t}{\tau}\right)\frac{d\tau}{\tau }\\
\phantom{f(t)}{}
= -\sum_{k\ge 0}\res\left\{g(\tau )\theta_p\!\left(\frac{t}{\tau}\right)\!\frac{1}{\tau };\tau \!=\!-p^{-2\nu-2-k}\right\}\!
 -\sum_{k\ge 0}\res\left\{g(\tau )\theta_p\!\left(\frac{t}{\tau}\right)\!\frac{1}{\tau };\tau \!=\! -p^{-2-k}\right\},
\end{gather*}
where $0<r<r_0$. To calculate the residue, we use the following lemma.
\begin{lemma}\label{lems}For any $k\in\mathbb{N}$, $\lambda\in\mathbb{C}^*$, we have
\begin{enumerate}\itemsep=0pt
\item[$1.$] $\res\left\{\dfrac{1}{\left(\tau /\lambda ;p\right)_\infty}\dfrac{1}{\tau }:\tau =\lambda p^{-k}\right\}
=\dfrac{(-1)^{k+1}p^{\frac{k(k+1)}{2}}}{(p;p)_k (p;p)_\infty}$,
\item[$2.$] $\dfrac{1}{(\lambda p^{-k};p)_\infty}=\dfrac{(-\lambda )^{-k}p^{\frac{k(k+1)}{2}}}{(\lambda ;p)_\infty \left(p/\lambda ;p\right)_k},\qquad \lambda \not \in p^{\mathbb{Z}}$.
\end{enumerate}
\end{lemma}

Summing up all of the residues, we obtain the convergent series $f(t)$ as follows
\begin{gather*}
f(t) =\frac{\theta
_p\left(-p^{2\nu +2}t\right)}{(p^{-2\nu},p;p)_\infty}
{}_1\varphi_1\left(0,p^{1+2\nu};p,x\right)
 +\frac{\theta _p\left(-p^2t\right)}{(p^{2\nu},p;p)_\infty}
{}_1\varphi_1\left(0,p^{1-2\nu};p,p^{-2\nu}x\right),
\end{gather*}
where $xt=1$.
Therefore, we acquire the connection formula for $z(t)=\mathcal{E}(t)f(t)$.

\section{The limit of the connection formula}\label{section3}

In this section, we show that the limit $p\to 1^-$ of the connection formula gives a connection formula of the Bessel function. At f\/irst, we assume that $0<p<1$ and $0<\sqrt{p}<1$. For the Bessel function, we set the Hankel function of the f\/irst and the second kind $H_\nu^{(1)}(z)$ and $H_\nu^{(2)}(z)$.
\begin{definition}The Hankel function of the f\/irst kind is given by
\[
H_\nu^{(1)}(z):=\frac{\Gamma \left(\frac{1}{2}-\nu\right)}{\pi i\sqrt{\pi}}\left(\frac{z}{2}\right)^\nu\int_{1+\infty i}^{(1+)}e^{izt}\left(t^2-1\right)^{\nu -\frac{1}{2}}dt,\qquad -\pi <\arg z<2\pi .
\]
The Hankel function of the second kind is def\/ined by
\[H_\nu^{(2)}(z):=\frac{\Gamma \left(\frac{1}{2}-\nu\right)}{\pi i\sqrt{\pi}}\left(\frac{z}{2}\right)^\nu\int_{-1+\infty i}^{(-1-)}e^{izt}\left(t^2-1\right)^{\nu -\frac{1}{2}}dt,\qquad -2\pi <\arg z<\pi .\]
\end{definition}
The contour for $H_\nu^{(1)}(z)$ is a path starting from $t=+1+\infty i$, rounding the circle around $t=1$ counterclockwise, and going back to $t=+1+\infty i$. Moreover, the contour for $H_\nu^{(2)}(z)$ is a~path starting from $t=-1+\infty i$, rounding the circle around $t=1$ clockwise, and going back to $t=-1+\infty i$.

The Hankel functions can be written by $J_\nu (z)$:
\begin{gather}
H_\nu^{(1)}(z) =\frac{ie^{-\nu\pi i}}{\sin\nu\pi}\left\{J_\nu (z)-e^{\nu\pi i}J_{-\nu}(z)\right\},\label{ha1}\\
H_\nu^{(2)}(z) =-\frac{ie^{\nu\pi i}}{\sin\nu\pi}\left\{J_\nu (z)-e^{-\nu\pi i}J_{-\nu}(z)\right\}.\label{ha2}
\end{gather}
The Hankel functions have asymptotic expansions around $z=0$ \cite{Ol}:
\begin{gather*}
H_\nu^{(1)}(z) \sim \left(\frac{2}{\pi z}\right)^{\frac{1}{2}}e^{i\zeta}\sum_{s\ge 0}i^s\frac{A_s(\nu )}{z^s},\qquad -\pi +\delta \leq \arg z\leq 2\pi -\delta ,\\
H_\nu^{(2)}(z) \sim \left(\frac{2}{\pi z}\right)^{\frac{1}{2}}e^{-i\zeta}\sum_{s\ge 0}(-i)^s\frac{A_s(\nu )}{z^s},\qquad -2\pi +\delta \leq \arg z\leq \pi -\delta ,
\end{gather*}
as $z\to \infty$. Here, $\delta $ is an any small constant,
\[A_s(\nu )=\frac{(4\nu^2-1^2)(4\nu^2-3^2)\cdots \left\{4\nu^2-(2s-1)^2\right\}}{s!8^s}\]
and
\[\zeta =z-\frac{1}{2}\nu\pi -\frac{1}{4}\pi .\]
In this sense, (\ref{ha1}) and (\ref{ha2}) considered as connection formula of the Bessel equation.

\subsection{Limit of the connection formula}\label{section3.1}

We rewrite the connection formula in Theorem~\ref{theorem1} in order to take a limit $p\to 1^-$. We set new functions $h_\nu (t;p)$ and $J_\nu^{\pm}(x;p)$. We set $h_\nu (t;p):=(p^{1/2},p^{1/2};p)_\infty z(t)$.
For any $x\in\mathbb{C}^*\setminus [-\lambda ;p]$ and $\lambda\in\mathbb{C}^*$, $J_{\nu ,\lambda}^{+}(x;p)$ is
\[J_{\nu ,\lambda}^{+}(x;p):=\frac{(p^{\nu +1};p)_\infty}{(p;p)_\infty}\frac{\theta_p\left(\frac{\lambda p^{\nu}}{x}\right)}{\theta_p\left(\frac{\lambda}{x}\right)}{}_1\varphi_1\left(0;p^{1+2\nu};p,x\right).\]
Similarly, $J_{\nu ,\lambda}^{-}(x;p)$ is
\[J_{\nu ,\lambda}^{-}(x;p):=\frac{(p^{\nu +1};p)_\infty}{(p;p)_\infty}\frac{\theta_p\left(\frac{\lambda p^{\nu}}{x}\right)}{\theta_p\left(\frac{\lambda}{x}\right)}{}_1\varphi_1\left(0;p^{1+2\nu};p,p^{2\nu }x\right).\]
We remark that the function $\theta_p(\lambda p^\nu /x)/\theta_p(\lambda /x)$ satisf\/ies the following $q$-dif\/ference equation
\[
u(px)=p^\nu u(x),
\]
which is also satisf\/ied by the function $u(x)=x^\nu$. We remark that the pair
$(J_{\nu,\lambda}^{+}(x;p),J_{-\nu,\lambda}^{-}(x;p))$ gives a fundamental system of solutions of equation (\ref{b3}) if $\nu\not\in\mathbb{Z}$. We set the function $C^+_\nu (\lambda ,t;p)$ and $C^-_\nu (\lambda ,t;p)$ as follow:
\begin{definition}For any $\lambda\in\mathbb{C}^*$, $C^+_\nu (\lambda ,t;p)$ is
\[C^+_\nu (\lambda ,t;p):=
\frac{(p^{\frac{1}{2}},p^{\frac{1}{2}};p)_\infty}{(p^{\nu +1},p^{-2\nu};p)_\infty}\frac{\theta_p(-p^{2\nu +2}t)}{\theta_p(-p^{\nu +2}t)}\frac{\theta_p(\lambda t)}{\theta_p(\lambda p^\nu t)}.\]
Similarly, the function $C^-_\nu (\lambda ,t;p)$ is
\[C^-_\nu (\lambda ,t;p):=
\frac{(p^{\frac{1}{2}},p^{\frac{1}{2}};p)_\infty}{(p^{-\nu +1},p^{2\nu};p)_\infty}\frac{\theta_p(-p^2t)}{\theta_p(-p^{\nu +2}t)}\frac{\theta_p(\lambda t)}{\theta_p(\lambda p^{-\nu }t)}.\]
\end{definition}
Then, $C^+_\nu (\lambda ,t;p)$ and $C^-_\nu (\lambda ,t;p)$ are single valued as a function of $t$. The function $C^+_\nu (\lambda ,t;p)$ and $C^-_\nu (\lambda ,t;p)$ are the $p$-elliptic functions.
By using these new functions, our connection formula is rewritten by
\[h_\nu\left(\frac{1}{x};p\right)=C_{\nu}^+\left(\lambda ,\frac{1}{x};p\right)J_{\nu}^{+}(x;p)
+C_\nu^-\left(\lambda ,\frac{1}{x};p\right)J_{-\nu ,\lambda}^{-}(x;p).\]
\begin{theorem}\label{limitbessel}For any $x\in \mathbb{C}^*\setminus (-\infty ,0]$  where $\arg x\in (-\pi ,\pi )$, we have
\[\lim_{p\to 1^-}h_\nu \left(\frac{1}{(1-p)^2x};p\right)
=-ie^{-\nu\pi i}H_{2\nu}^{(2)}(2\sqrt{x}).\]
Here, $H_{2\nu}^{(2)}(\cdot )$ is the Hankel function of the second kind.
\end{theorem}
The aim of this section is to give a proof of the theorem above.

By the def\/inition, $h_\nu \left(1/\{(1-p)^2x\};p\right)$ can be described as follows
\begin{gather}
 h_\nu\left(\frac{1}{(1-p)^2x};p\right)
 =\left\{\frac{(p^{\frac{1}{2}},p^{\frac{1}{2}};p)_\infty}{(p^{-2\nu},p;p)_\infty}(1-p)^{2\nu}\right\}
\left\{\frac{\theta_p\left(-\frac{p^{2\nu +2}}{x(1-p)^2}\right)}{\theta_p\left(-\frac{p^{\nu +2}}{x(1-p)^2}\right)}(1-p)^{-2\nu}\right\} \notag.\\
 \hphantom{h_\nu\left(\frac{1}{(1-p)^2x};p\right) =}{}
  \times\left\{{}_1\varphi_1\left(0;p^{1+2\nu};p,(1-p)^2x\right)\right\}\notag\\
\hphantom{h_\nu\left(\frac{1}{(1-p)^2x};p\right) =}{}  +
\left\{\frac{(p^{\frac{1}{2}},p^{\frac{1}{2}};p)_\infty}{(p^{2\nu},p;p)_\infty}(1-p)^{-2\nu}\right\}
\left\{\frac{\theta_p\left(-\frac{p^2}{x(1-p)^2}\right)}{\theta_p\left(-\frac{p^{\nu +2}}{x(1-p)^2}\right)}(1-p)^{2\nu}\right\} \notag\\
\hphantom{h_\nu\left(\frac{1}{(1-p)^2x};p\right) =}{} \times
\left\{{}_1\varphi_1\left(0;p^{1-2\nu};p,p^{-2\nu}(1-p)^2x\right)\right\}.\label{limit}
\end{gather}
We consider the limit of each part $\{\cdot \}$.

\begin{lemma}\label{limgam}For any $\nu \in\mathbb{C}^*\setminus \mathbb{Z}$, we have
\[\lim_{p\to 1^-}
\frac{(p^{\frac{1}{2}},p^{\frac{1}{2}};p)_\infty}{(p^{-2\nu},p;p)_\infty}(1-p)^{2\nu}
=-\frac{1}{\sin (2\nu\pi )\Gamma (2\nu +1)}.
\]
\end{lemma}
\begin{proof}We can check out as follows
\begin{gather*}
\frac{(p^{\frac{1}{2}},p^{\frac{1}{2}};p)_\infty}{(p^{-2\nu},p;p)_\infty}(1-p)^{2\nu} =\frac{\frac{(p;p)_\infty}{(p^{-2\nu};p)_\infty}(1-p)^{1+2\nu}}
{\left\{\frac{(p;p)_\infty}{(p^{\frac{1}{2}};p)_\infty}(1-p)^{\frac{1}{2}}\right\}
\left\{\frac{(p;p)_\infty}{(p^{\frac{1}{2}};p)_\infty}(1-p)^{\frac{1}{2}}\right\}
} =\frac{\Gamma_p(-2\nu )}{\Gamma_p\left(\frac{1}{2}\right)\Gamma_p\left(\frac{1}{2}\right)}.
\end{gather*}
Here, $\Gamma_q(\cdot )$ is Jackson's $q$-gamma function which is def\/ined by
\[\Gamma_q(x):=\frac{(q;q)_\infty}{(q^x;q)_\infty}(1-q)^{1-x},\qquad 0<q<1.\]
This function satisf\/ies $\lim\limits_{q\to 1^-}\Gamma_q(x)=\Gamma (x)$~\cite{GR}. Therefore,
\[\lim_{p\to 1^-}
\frac{(p^{\frac{1}{2}},p^{\frac{1}{2}};p)_\infty}{(p^{-2\nu},p;p)_\infty}(1-p)^{2\nu}
=\frac{\Gamma (-2\nu )}{\Gamma \left(\frac{1}{2}\right)\Gamma \left(\frac{1}{2}\right)}.\]
By Euler's ref\/lection formula of the gamma function, we get
\[\frac{\Gamma (-2\nu )}{\Gamma\left(\frac{1}{2}\right)\Gamma\left(\frac{1}{2}\right)}=
-\frac{1}{\sin (2\nu\pi )\Gamma (2\nu +1)}.\]
Therefore, we get the conclusion.
\end{proof}
If we replace $\nu$ by $-\nu$, we get the limit
\[
\lim_{p\to 1^-}
\frac{(p^{\frac{1}{2}},p^{\frac{1}{2}};p)_\infty}{(p^{2\nu},p;p)_\infty}(1-p)^{-2\nu}=\frac{1}{\sin (2\nu\pi )\Gamma (1-2\nu )}.
\]
In \cite{Z1}, the following proposition can be found:

\begin{proposition}\label{propo1}
 For any $x\in\mathbb{C}^*$ $(-\pi <\arg x<\pi )$, we have
\[\lim_{p\to 1^-}\frac{\theta_p\left(\frac{p^{\nu_1}}{(1-p^2)x}\right)}{\theta_p\left(\frac{p^{\nu_2}}{(1-p^2)x}\right)}\left(1-p^2\right)^{\nu_2-\nu_1}=x^{\nu_1-\nu_2},\]and
\[
\lim_{p\to 1^-}\frac{\theta_p\left(-\frac{p^{\nu_1}}{(1-p^2)x}\right)}{\theta_p\left(-\frac{p^{\nu_2}}{(1-p^2)x}\right)}\left(1-p^2\right)^{\nu_2-\nu_1}=(-x)^{\nu_1-\nu_2}.
\]
\end{proposition}

\begin{lemma}For any $x\in\mathbb{C}^*$ $(-\pi <\arg x\leq\pi )$ and fixed constant~$K$, we have
\[\theta_p(-\sqrt{p})\theta_p\left(-\frac{K}{x}\right)
=\theta_{\sqrt{p}}\left(\sqrt{\frac{K}{x}}\right)\theta_{\sqrt{p}}\left(-\sqrt{\frac{K}{x}}\right).\]
\end{lemma}

\begin{proof}From Jacobi's triple product identity and $(a^2;q^2)_n=(a,-a;q)_n$, we obtain
\[\frac{(\sqrt{p};\sqrt{p})_\infty}{(-\sqrt{p};\sqrt{p})_\infty}\theta_p\left(-\frac{K}{x}\right)
=\theta_{\sqrt{p}}\left(\sqrt{\frac{K}{x}}\right)\theta_{\sqrt{p}}\left(-\sqrt{\frac{K}{x}}\right).\]
We remark that $(\sqrt{p};\sqrt{p})_\infty /(-\sqrt{p};\sqrt{p})_\infty$ can be rewritten as follows \cite{GR}:
\[\frac{(\sqrt{p};\sqrt{p})_\infty}{(-\sqrt{p};\sqrt{p})_\infty}
=\sum_{n\in\mathbb{Z}}(-1)^n(\sqrt{p})^{n^2}=\theta_p(-\sqrt{p}).\]
We obtain the conclusion.
\end{proof}
Therefore, we obtain the following relation.
\begin{corollary}\label{kei2}For any $x\in\mathbb{C}^*$ $(-\pi <\arg x\leq\pi )$, we have
\begin{gather}
\frac{\theta_p\left(p^{2\nu +2}\frac{-1}{(1-p)^2x}\right)}{\theta_p\left(p^{\nu +2}\frac{-1}{(1-p)^2x}\right)}
=\frac{\theta_{\sqrt{p}}\left(p^{\nu +1}\frac{1}{(1-p)\sqrt{x}}\right) \theta_{\sqrt{p}}\left(p^{\nu +1}\frac{-1}{(1-p)\sqrt{x}}\right)}
{\theta_{\sqrt{p}}\left(p^{\frac{\nu}{2}+1}\frac{1}{(1-p)\sqrt{x}}\right)\theta_{\sqrt{p}}
\left(p^{\frac{\nu}{2}+1}\frac{-1}{(1-p)\sqrt{x}}\right)}\label{th1}
\end{gather}
and
\begin{gather}
\frac{\theta_p\left(p^2\frac{-1}{(1-p)^2x}\right)}{\theta_p\left(p^{\nu +2}\frac{-1}{(1-p)^2x}\right)}
=\frac{\theta_{\sqrt{p}}\left(p\frac{1}{(1-p)\sqrt{x}}\right)\theta_{\sqrt{p}}\left(p\frac{-1}{(1-p)\sqrt{x}}\right)}
{\theta_{\sqrt{p}}\left(p^{\frac{\nu}{2}+1}\frac{1}{(1-p)\sqrt{x}}\right)\theta_{\sqrt{p}}\left(p^{\frac{\nu}{2}+1}\frac{-1}{(1-p)\sqrt{x}}\right)}.\label{th2}
\end{gather}
\end{corollary}
We consider the limit $p\to 1^-$ (i.e., $\sqrt{p}\to 1^-$) of \eqref{th1} and \eqref{th2}.
\begin{lemma}For any $x\in\mathbb{C}^*\setminus (-\infty ,0]$ $(-\pi <\arg x\leq\pi )$, we have
\label{limthe}
\begin{enumerate}\itemsep=0pt
\item[$1.$] $\displaystyle{\lim_{p\to 1^-}\frac{\theta_p\left(-\frac{p^{2\nu +2}}{x(1-p)^2}\right)}
{\theta_p\left(-\frac{p^{\nu +2}}{x(1-p)^2}\right)}(1-p)^{-2\nu}=e^{\nu \pi i}x^\nu}$ \
and
\item[$2.$] \label{limtheta}$\displaystyle{\lim_{p\to 1^-}\frac{\theta_p\left(-\frac{p^{2}}{x(1-p)^2}\right)}
{\theta_p\left(-\frac{p^{\nu +2}}{x(1-p)^2}\right)}(1-p)^{2\nu}=e^{-\nu \pi i}x^{-\nu}}$.
\end{enumerate}
\end{lemma}

\begin{proof}
Combining Proposition~\ref{propo1} and Corollary~\ref{kei2}, we consider the limit $\sqrt{p}\to 1^-$ as follows:
\begin{gather*}
 \frac{\theta_p\left(p^{2\nu +2}\frac{-1}{(1-p)^2x}\right)}{\theta_p\left(p^{\nu +2}\frac{-1}{(1-p)^2x}\right)}(1-p)^{-2\nu}
=\frac{\theta_{\sqrt{p}}\left(p^{\nu +1}\frac{1}{(1-p)\sqrt{x}}\right) \theta_{\sqrt{p}}\left(p^{\nu +1}\frac{-1}{(1-p)\sqrt{x}}\right)}
{\theta_{\sqrt{p}}\left(p^{\frac{\nu}{2}+1}\frac{1}{(1-p)\sqrt{x}}\right)\theta_{\sqrt{p}}\left(p^{\frac{\nu}{2}+1}\frac{-1}{(1-p)\sqrt{x}}\right)}
(1-p)^{-2\nu}\\
 \hphantom{\frac{\theta_p\left(p^{2\nu +2}\frac{-1}{(1-p)^2x}\right)}{\theta_p\left(p^{\nu +2}\frac{-1}{(1-p)^2x}\right)}(1-p)^{-2\nu}}{}
 =\left\{\frac{\theta_{\sqrt{p}}\left((\sqrt{p})^{2\nu +2}\frac{1}{(1-(\sqrt{p})^2)\sqrt{x}}\right) }{\theta_{\sqrt{p}}\left((\sqrt{p})^{\nu +2}\frac{1}{(1-\{\sqrt{p})^2)\sqrt{x}}\right)}\left\{1-(\sqrt{p})^2\right\}^{-\nu}\right\}\\
\hphantom{\frac{\theta_p\left(p^{2\nu +2}\frac{-1}{(1-p)^2x}\right)}{\theta_p\left(p^{\nu +2}\frac{-1}{(1-p)^2x}\right)}(1-p)^{-2\nu}=}{}
\times
\left\{\frac{\theta_{\sqrt{p}}\left(-(\sqrt{p})^{2\nu +2}\frac{1}{(1-(\sqrt{p})^2)\sqrt{x}}\right) }{\theta_{\sqrt{p}}\left(-(\sqrt{p})^{\nu +2}\frac{1}{(1-\{\sqrt{p})^2)\sqrt{x}}\right)}\left\{1-(\sqrt{p})^2\right\}^{-\nu}\right\}\\
\hphantom{\frac{\theta_p\left(p^{2\nu +2}\frac{-1}{(1-p)^2x}\right)}{\theta_p\left(p^{\nu +2}\frac{-1}{(1-p)^2x}\right)}(1-p)^{-2\nu}}{}
 \to (\sqrt{x})^\nu\cdot (-\sqrt{x})^\nu=(-x)^\nu
=e^{\nu\pi i}x^\nu , \qquad \sqrt{p}\to 1^-.
\end{gather*}
Similarly, we can prove the latter one. We obtain the conclusion.
\end{proof}

We consider the last part.
\begin{lemma}\label{limbas}For any $x\in\mathbb{C}^*$, we have
\[
\lim_{p\to 1^-}{}_1\varphi_1\left(0;p^{1+2\nu};p,(1-p)^2x\right)={}_0F_1\left(-,1+2\nu ;-x\right)\]
and
\[
\lim_{p\to 1^-}{}_1\varphi_1\left(0;p^{1-2\nu};p,p^{-2\nu }(1-p)^2x\right)={}_0F_1\left(-,1-2\nu ;-x\right).
\]
\end{lemma}
\begin{proof}
We check each of the term of
\[{}_1\varphi_1\left(0;p^{1+2\nu};p,(1-p)^2x\right)
=\sum_{n\ge 0}\frac{1}{(p^{1+2\nu},p;p)_n}(-1)^np^{\frac{n(n-1)}{2}}\left\{(1-p)^2x\right\}^n.\]
For any $n\ge 0$,
\begin{gather*}
\frac{1}{(p^{1+2\nu},p;p)_n}(-1)^np^{\frac{n(n-1)}{2}}
\left\{(1-p)^2x\right\}^n \\
\qquad {}=\frac{(1-p)^n(1-p)^n}{(p^{1+2\nu};p)_n(p;p)_n}p^{\frac{n(n-1)}{2}}\left(-x\right)^n
 \to \frac{1}{(1+2\nu )_n\cdot n!}\left(-x\right)^n, \qquad p\to 1^-.
\end{gather*}
Summing up all terms, we get
\[\sum_{n\ge 0}\frac{1}{(1+2\nu )_n\cdot n!}\left(-x\right)^n
={}_0F_1\left(-,1+2\nu ;-x\right).\]
Therefore, we obtain the conclusion. Similarly, we can prove the latter.
\end{proof}

We give the proof of Theorem~\ref{limitbessel}.

\begin{proof}Apply Lemma \ref{limgam}, Lemma \ref{limthe} and Lemma \ref{limbas} to (\ref{limit}), we obtain
\begin{gather*}
h_\nu \left(\frac{1}{(1-p)^2x};p\right)\to  \left\{-\frac{1}{\sin (2\nu\pi )\Gamma (1+2\nu )}\right\}e^{\nu\pi i}x^{\nu}{}_0F_1\left(-,1+2\nu ;-x\right)\\
\hphantom{h_\nu \left(\frac{1}{(1-p)^2x};p\right)\to}{}
+\left\{\frac{1}{\sin (2\nu\pi )\Gamma (1-2\nu )}\right\}e^{-\nu\pi i}x^{-\nu}{}_0F_1\left(-,1-2\nu ;-x\right)\\
\hphantom{h_\nu \left(\frac{1}{(1-p)^2x};p\right)}{}
=\frac{-e^{\nu\pi i}J_{2\nu}\left(2\sqrt{x}\right)+
e^{-\nu\pi i}J_{-2\nu}\left(2\sqrt{x}\right)}{\sin (2\nu\pi )}\\
\hphantom{h_\nu \left(\frac{1}{(1-p)^2x};p\right)}{}
=\frac{e^{-\nu\pi i}}{i}H_{2\nu}^{(2)}\left(2\sqrt{x}\right)
,\qquad p\to 1^-.
\end{gather*}
Therefore, we acquire the conclusion.
\end{proof}

\subsection*{Acknowledgements}
The author would like to express his deepest gratitude to Professor Yousuke Ohyama for many valuable comments. The author also expresses his thanks to Professor Lucia Di Vizio for fruitful discussions when she was invited to the University of Tokyo in the winter 2011. The author would like to give thanks to the referee for some useful comments.

\pdfbookmark[1]{References}{ref}
\LastPageEnding

\end{document}